\documentclass[12pt]{article}
\usepackage{amssymb}
\usepackage{amsthm}
\usepackage{fancyhdr}
\usepackage{amscd}
\usepackage{extarrows}
\usepackage{graphicx}
\usepackage{fancyhdr}
\usepackage{epsfig}
\usepackage{amsfonts}
\usepackage{mathrsfs}
\usepackage{picinpar} \usepackage{amsmath} \usepackage[all]{xy}
\usepackage[colorlinks=true]{hyperref}
\hypersetup{urlcolor=blue, citecolor=red}
\theoremstyle{plain}
\newtheorem{theorem}{Theorem}[section]
\newtheorem{corollary}[theorem]{Corollary}
\newtheorem{lemma}[theorem]{Lemma}

\newtheorem{condition}[theorem]{Condition}
\newtheorem{proposition}[theorem]{Proposition}
\theoremstyle{definition}

\newtheorem{definition}[theorem]{Definition}

\newtheorem{remark}[theorem]{Remark}

\newcommand{\Add}{\mbox{\rm Add}}

\newcommand{\xt}{\mbox{\rm xt}}
\newcommand{\for}{\mbox{\rm for}}
\newcommand{\any}{\mbox{\rm any}}

\newcommand{\op}{\mbox{\rm op}}
\newcommand{\Pres}{\mbox{\rm Pres}}

\begin{document}
\title{$\xi$-tilting objects in extriangulated categories\thanks{Supported by the National Natural Science Foundation of China (Grant No. 11771212, 11801004)
and a project funded by the Priority Academic Program Development of Jiangsu Higher Education Institutions.}}
\smallskip
\author{\small Yuxia Mei and Jiaqun Wei\\%
\small E-mail addresses: meiyuxia2010@163.com, weijiaqun@njnu.edu.cn\\
\small Institute of Mathematics, School of Mathematics Sciences\\
\small Nanjing Normal University, 
Nanjing 210023, P.R.China\\
\small 
}

\date{}
\maketitle


\noindent{\footnotesize {{\textbf{Abstract}}: Extriangulated categories were introduced by Nakaoka and Palu by extracting the similarities between exact categories and triangulated categories. In this article, we introduce and study the notion of $\xi$-tilting object in an extriangulated category, where $\xi$ is a proper class of $\mathbb{E}$-triangles. This extends results in \cite{YGH_2020}.
\vspace{0.2cm}

\noindent{\footnotesize {{\textbf{MSC2010}}: 18A40 18E10 18E30}

\noindent{\footnotesize {{\textbf{Keywords}}: Extriangulated category; proper class; $\xi$-projective; $\xi$-tilting object }

\section{Introduction}
The classical tilting theory was introduced in the context of finite generated modules over artin algebras by Auslander et al. \cite{MA_1979}, Brenner and Butler \cite{SB_1980}, Happel and Ringel \cite{DH_1982} and so on, and since then it has played a central role in the development of the representation theory of artin algebras.

Triangulated categories were introduced by Verdier \cite{JLV_1997}. Keller \cite{BK_2007} introduced the notion of tilting objects in an algebra triangulated category. The homological theory of triangulated categories was developed by Beligiannis \cite{AB_2000}. It parallels the homological theory in an exact category in the sense of Quillen. Using the proper class $\xi$ of triangles in a triangulated category $\mathcal{C}$, Beligiannis defined the notions of $\xi$-projective objects, $\xi$-injective objects, $\xi$-projective dimension, and so on. Recently, Y.G. Hu et al. \cite{YGH_2020} introduced the notion of $\xi$-tilting object in a triangulated category by means of the homological theory in the triangulated category.

Triangulated categories and exact categories are two fundamental structures in mathematics. They are also important tools in many mathematics branches. It is well known that these two kinds of categories have some similarities: while exact categories admit short exact sequences, triangulated categories admit triangles. By extracting the similarities between triangulated categories and exact categories, Nakaoka and Palu \cite{NP_2019} introduced the notion of extriangulated categories. Such category is a triplet  $(\mathcal{C},\mathbb{E},\mathfrak{s})$, where $\mathcal{C}$ is an additive category, $\mathbb{E}: \mathcal{C}^{op}\times \mathcal{C}\rightarrow \mathcal{A}b$ is a biadditive functor and $\mathfrak{s}$ assigns to each $\delta \in \mathbb{E}(C,A)$ a class of 3-terms sequences with end terms $A$ and $C$ such that certain axioms hold. Using the proper class $\xi$ of $\mathbb{E}$-triangles in an extriangulated category $\mathcal{C}$, J.S. Hu \cite{JSH1_2020} defined the notions of $\xi$-projective objects, $\xi$-injective objects, $\xi$-projective dimension, and so on.

The aim of this paper is to extend the results of Y.G. Hu et al. \cite{YGH_2020} to extriangulated categories. We introduce the notion of $\xi$-tilting object in a triangulated category and study their properties. More precisely, we obtain the $\mathbb{E}$-triangle versions of the Bazzoni characterization in extriangulated category.

This paper is organized as follows. In section 2, we recall the definition of extriangulated categories and some basic properties which are needed in the paper. In section 3, we extend Y.G. Hu's notion of $\xi$-tilting object in a triangulated category\cite{YGH_2020} to $\xi$-tilting object in an extriangulated category and we state and prove our main result. We get the Bazzoni characterization of $\xi$-tilting object in extriangulated category as follows. For unexplained notions in the following theorem, we refer to Section $2$.

$\mathbf{Theorem}$ \ref{T34} Let $T$ be an object in $\mathcal{C}$. Then the following statements are equivalent.

(1) $T$ is an $\xi$-tilting object.

(2) $T^{\perp} = \Pres_{\mathcal{P}(\xi)}^{1}(\Add T)$.

\section{Preliminaries}

Throughout this paper, let $\mathcal{C}$ be an additive category. For the category $\mathcal{C}$, if $A,B \in \mathcal{C}$, then we denote the set of morphisms $A \rightarrow B$ in $\mathcal{C}$ by $\mathcal{C}(A,B)$. We denote the identity morphism of an object $C \in \mathcal{C}$ by $1=1_{C}$. If $f \in \mathcal{C}(A,B)$, $g \in \mathcal{C}(B,C)$, we denote composition of $f$ and $g$ by $gf$. By a subcategory $\mathcal{D}$ of $\mathcal{C}$, we always mean that $\mathcal{D}$ is an additive full subcategory which is closed under isomorphisms, direct sums and direct summands.

We recall some basics on extriangulated categories from \cite{NP_2019, JSH1_2020, JSH2_2020, JSH3_2020}.

Suppose $\mathcal{C}$ is equipped with a biadditive functor $\mathbb{E}:\mathcal{C}^{op}\times \mathcal{C}\rightarrow \mathcal{   A}b$, where $\mathcal{A}b$ is the category of abelian groups. For any pair of objects $A,C\in \mathcal{C}$, an element $\delta \in \mathbb{E}(C,A)$ is called an $\mathbb{E}$-$extension$. Thus formally, an $\mathbb{E}$-extension is a triplet $(A,\delta,C)$. Let $(A,\delta,C)$ be an $\mathbb{E}$-extension. Since $\mathbb{E}$ is a bifunctor, for any $a \in \mathcal{C}(A,A')$ and $c \in \mathcal{C}(C',C)$, we have $\mathbb{E}$-extensions
$$\mathbb{E}(C,a)(\delta)\in \mathbb{E}(C,A')~~~~\text{and}~~~~\mathbb{E}(c,A)(\delta)\in \mathbb{E}(C',A).$$
They are abbreviated to $a_{\ast}\delta$ and $c^{\ast}\delta$ respectively. In this terminology, we have $$\mathbb{E}(c,a)(\delta)=c^{\ast}a_{\ast}\delta=a_{\ast}c^{\ast}\delta$$
in $\mathbb{E}(C',A')$. For any $A,C \in \mathcal{C}$, the zero element $0 \in \mathbb{E}(C,A)$ is called the $split$ $\mathbb{E}$-$extension$.

\begin{definition} \cite[Definition 2.3]{NP_2019} Let $(A,\delta,C)$, $(A',\delta',C')$ be any pair of $\mathbb{E}$-extensions. A $morphism$
$$(a,c):(A,\delta ,C)\rightarrow (A',\delta',C')$$
of $\mathbb{E}$-extensions is a pair of morphisms $a \in \mathcal{C}(A,A')$ and $c \in \mathcal{C}(C,C')$ in $\mathcal{C}$, satisfying the equality $$a_{\ast}\delta=c^{\ast}\delta'.$$
Simply we denote it as $(a,c):\delta \rightarrow \delta'$.
\end{definition}

\begin{definition} \cite[Definition 2.6]{NP_2019} Let $\delta=(A,\delta,C),~\delta'=(A',\delta',C')$ be any pair of $\mathbb{E}$-extensions. Let $$C\stackrel{l_{C}}\rightarrow C \oplus C' \stackrel{l_{C'}}\leftarrow C'$$ and $$A\stackrel{p_{A}}\leftarrow A\oplus A'\stackrel{p_{A'}}\rightarrow A'$$
be coproduct and product in $\mathcal{C}$, respectively. Remark that, by the additivity of $\mathbb{E}$, we have a natural isomorphism $$\mathbb{E}(C\oplus C',A\oplus A')\simeq \mathbb{E}(C,A)\oplus \mathbb{E}(C,A')\oplus \mathbb{E}(C',A)\oplus \mathbb{E}(C',A'). $$

Let $\delta \oplus \delta'\in \mathbb{E}(C\oplus C',A\oplus A') $ be the element corresponding to $(\delta,0,0,\delta')$ through this isomorphism. This is the unique element which satisfies $$\mathbb{E}(l_{C},p_{A})(\delta \oplus \delta')=\delta,~~~~\mathbb{E}(l_{C},p_{A'})(\delta \oplus \delta')=0,$$
$$\mathbb{E}(l_{C'},p_{A})(\delta \oplus \delta')=0,~~~~\mathbb{E}(l_{C'},p_{A'})(\delta \oplus \delta')=\delta'.$$
\end{definition}

\begin{definition} \cite[Definition 2.7]{NP_2019} Let $A,C \in \mathcal{C}$ be any pair of objects. Sequences of morphisms in $\mathcal{C}$ $$A\stackrel{x}\rightarrow B\stackrel{y}\rightarrow C~~~~\text{and}~~~~A\stackrel{x'}\rightarrow B'\stackrel{y'}\rightarrow C$$
are said to be $equivalent$ if there exists an isomorphism $b\in \mathcal{C}(B,B')$ which makes the following diagram commutative.

$$\xymatrix{
 A\ar[r]^{x}\ar@{=}[d]&B\ar[r]^{y}\ar[d]_{\simeq}^{b}& C\ar@{=}[d] \\
 A\ar[r]_{x'}& B'\ar[r]_{y'}      &C
  }
$$
We denote the equivalence class of $A\stackrel{x}\rightarrow B\stackrel{y}\rightarrow C$ by $[A\stackrel{x}\rightarrow B\stackrel{y}\rightarrow C]$. For any $A,C\in \mathcal{C}$, we denote as $0=[A\stackrel{\left[\begin{smallmatrix} \  1_{A} \\ 0 \end{smallmatrix}\right]}\rightarrow A\oplus C\stackrel{\left[\begin{smallmatrix} \  0&1_{C} \end{smallmatrix}\right]}\rightarrow C]$.

\end{definition}

\begin{definition} \cite[Definition 2.9 and 2.10]{NP_2019} Let $\mathfrak{s}$ be a correspondence which associates an equivalence class $\mathfrak{s}(\delta)=[A\stackrel{x}\rightarrow B\stackrel{y}\rightarrow C]$ to any $\mathbb{E}$-extension $\delta \in \mathbb{E}(C,A)$. This $\mathfrak{s}$ is called a $realization$ of $\mathbb{E}$, if it satisfies the following condition:

$(\bullet)$ Let $\delta \in \mathbb{E}(C,A)$ and $\delta' \in \mathbb{E}(C',A')$ be any pair of $\mathbb{E}$-extensions, with $$\mathfrak{s}(\delta)=[A\stackrel{x}\rightarrow B\stackrel{y}\rightarrow C],~~~~\mathfrak{s}(\delta')=[A'\stackrel{x'}\rightarrow B'\stackrel{y'}\rightarrow C'].$$
Then, for any morphism $(a,c):\delta\rightarrow \delta'$, there exists $b \in \mathcal{C}(B,B')$ which makes the following diagram commutative.

$$\xymatrix{
 A\ar[r]^{x}\ar[d]^{a}&B\ar[r]^{y}\ar[d]^{b}& C\ar[d]^{c} \\
 A'\ar[r]_{x'}& B'\ar[r]_{y'}      &C'
  }
$$
In this case, we say that sequence $A\stackrel{x}\rightarrow B\stackrel{y}\rightarrow C$ realizes $\delta$, whenever it satisfies $\mathfrak{s}(\delta)=[A\stackrel{x}\rightarrow B\stackrel{y}\rightarrow C]$. Note that this condition does not depend on the choices of the representatives of the equivalence classes. In the above situation, we say that $(a,b,c)$ $realizes$ $(a,c)$.

A realization $\mathfrak{s}$ of $\mathbb{E}$ is called $additive$ if the following conditions are satisfied.

$(1)$ For any $A,C\in \mathcal{C}$, the split $\mathbb{E}$-extension $0\in \mathbb{E}(C,A)$ satisfies $\mathfrak{s}(0)=0.$

$(2)$ For any pair of $\mathbb{E}$-extensions $\delta \in \mathbb{E}(C,A)$ and $\delta' \in \mathbb{E}(C',A')$, we have $\mathfrak{s}(\delta \oplus \delta')=\mathfrak{s}(\delta)\oplus \mathfrak{s}(\delta')$.
\end{definition}

\begin{definition}\cite[Definition 2.12]{NP_2019} We call the triplet $(\mathcal{C},\mathbb{E},\mathfrak{s})$ an $externally$ $triangulated$ $category$ (or $extriangulated$ $category$ $\mathcal{C}$ for short) if it satisfies the following conditions:

(ET1) $\mathbb{E}:\mathcal{C}^{op}\times \mathcal{C}\rightarrow \mathcal{A}b$ is a biadditive functor.

(ET2) $\mathfrak{s}$ is an additive realization of $\mathbb{E}$.

(ET3) Let $\delta\in \mathbb{E}(C,A)$ and $\delta'\in \mathbb{E}(C',A')$ be any pair of $\mathbb{E}$-extensions, realized as
$$\mathfrak{s}(\delta)=[A\stackrel{x}\rightarrow B\stackrel{y}\rightarrow C],~~~~\mathfrak{s}(\delta')=[A'\stackrel{x'}\rightarrow B'\stackrel{y'}\rightarrow C'].$$
For any commutative square
$$\xymatrix{
 A\ar[r]^{x}\ar[d]^{a}&B\ar[r]^{y}\ar[d]^{b}& C \\
 A'\ar[r]_{x'}& B'\ar[r]_{y'}      &C'
  }
$$
in $\mathcal{C}$, there exists a morphism $(a,c):\delta \rightarrow \delta'$ which is realized by $(a,b,c)$.

(ET3)$^{\text{op}}$ Let $\delta\in \mathbb{E}(C,A)$ and $\delta'\in \mathbb{E}(C',A')$ be any pair of $\mathbb{E}$-extensions, realized as
$$\mathfrak{s}(\delta)=[A\stackrel{x}\rightarrow B\stackrel{y}\rightarrow C],~~~~\mathfrak{s}(\delta')=[A'\stackrel{x'}\rightarrow B'\stackrel{y'}\rightarrow C'].$$
For any commutative square
$$\xymatrix{
 A\ar[r]^{x}&B\ar[r]^{y}\ar[d]^{b}& C\ar[d]^{c} \\
 A'\ar[r]_{x'}& B'\ar[r]_{y'}      &C'
  }
$$
in $\mathcal{C}$, there exists a morphism $(a,c):\delta \rightarrow \delta'$ which is realized by $(a,b,c)$.

(ET4) Let $(A,\delta,D)$ and $(B,\delta',F)$ be $\mathbb{E}$-extensions realized by $$A\stackrel{f}\rightarrow B\stackrel{f'}\rightarrow C,~~~~B\stackrel{g}\rightarrow C\stackrel{g'}\rightarrow F$$ respectively. Then there exist an object $E \in \mathcal{C}$, commutative diagram
$$\xymatrix{
 A\ar@{=}[d]\ar[r]^{f}&B\ar[r]^{f'}\ar[d]^{g}& C\ar[d]^{d} \\
 A\ar[r]_{h}& C\ar[r]_{h'}\ar[d]^{g'}      &E\ar[d]^{e}\\
 &F\ar@{=}[r]&F
  }
$$
in $\mathcal{C}$, and an $\mathbb{E}$-extension $\delta'' \in \mathbb{E}(E,A)$ realized by $A\stackrel{h}\rightarrow C\stackrel{h'}\rightarrow E$, which satisfy the following compatibilities.

(i) $D\stackrel{d}\rightarrow E\stackrel{e}\rightarrow E$ realizes $f'_{\ast}\delta'$,

(ii) $d^{\ast}\delta''=\delta$,

(iii) $f_{\ast}\delta''=e^{\ast}\delta'$.

(ET4)$^{\text{op}}$ Let $(D,\delta,B)$ and $(F,\delta',C)$ be $\mathbb{E}$-extensions realized by $$D\stackrel{f'}\rightarrow A\stackrel{f}\rightarrow B~~~~\text{and}~~~~F\stackrel{g'}\rightarrow B\stackrel{g}\rightarrow C$$ respectively. Then there exists an object $E \in \mathcal{C}$, commutative diagram

$$\xymatrix{
 D\ar@{=}[d]\ar[r]^{d}&E\ar[r]^{e}\ar[d]^{h'}& F\ar[d]^{g'} \\
 D\ar[r]_{f'}& A\ar[r]_{f}\ar[d]^{h}      &B\ar[d]^{g}\\
 &C\ar@{=}[r]&C
  }
$$
in $\mathcal{C}$, and an $\mathbb{E}$-extension $\delta''\in \mathbb{E}(C,E)$ realized by $E\stackrel{h'}\rightarrow A\stackrel{h}\rightarrow C$, which satisfy the following compatibilities.

(i) $D\stackrel{D}\rightarrow E\stackrel{e}\rightarrow F$ realizes $g'^{\ast}\delta$,

(ii) $\delta'=e_{\ast}\delta''$,

(iii) $d_{\ast}\delta=g^{\ast}\delta''$.
\end{definition}

The following condition is analogous to the weak idempotent completeness in exact category (see \cite[Condition 5.8]{NP_2019}).

\begin{condition} (Condition (WIC))\label{W1}Consider the following conditions.

(1) Let $f \in \mathcal{C}(A,B)$, $g \in \mathcal{C}(B,C)$ be any composable pair of morphisms. If $gf$ is an inflation, then so is $f$.

(2) Let $f \in \mathcal{C}(A,B)$, $g \in \mathcal{C}(B,C)$ be any composable pair of morphisms. If $gf$ is an deflation, then so is $g$.

\end{condition}

\begin{lemma}\cite[Corollary 3.5]{NP_2019}\label{L3} Let $\mathcal{C}$ be an extriangulated category, and
$$\xymatrix{
 A\ar[r]^{x}\ar[d]^{a}&B\ar[r]^{y}\ar[d]^{b}& C\ar[d]^{c}\ar@{-->}[r]^{\delta}& \\
 A'\ar[r]_{x'}& B'\ar[r]_{y'}      &C'\ar@{-->}[r]^{\delta'}&
  }
$$
be any morphism of $\mathbb{E}$-triangles. Then the following are equivalent.

(1) $a$ factors through $x$;

(2) $a_{\ast}\delta=c^{\ast}\delta'=0$;

(3) $c$ factors through $y'$.

In particular, in this case $\delta=\delta'$ and $(a,b,c)=(1,1,1)$, we obtain $$x~is~a~section~\Leftrightarrow \delta~splits~\Leftrightarrow ~y~is~a~retraction.$$
The full subcategory consisting of the split $\mathbb{E}$-triangles will be denoted by $\Delta_{0}$.
\end{lemma}

The following concepts are quoted verbatim from \cite{JSH1_2020, JSH2_2020, JSH3_2020}. A class of $\mathbb{E}$-triangles $\xi$ is $closed$ $under$ $base~ change $ if for any $\mathbb{E}$-triangle
$$A\stackrel{x}\longrightarrow B\stackrel{y}\longrightarrow C\stackrel{\delta}\dashrightarrow$$ in $\xi$
and any morphism  $c:C' \rightarrow C$, then any $\mathbb{E}$-triangle
$A\stackrel{x'}\longrightarrow B'\stackrel{y'}\longrightarrow C'\stackrel{c^{*}\delta}\dashrightarrow$ belongs to $\xi$.

Dually, a class of $\mathbb{E}$-triangles $\xi$ is $closed~ under~ cobase~ change $ if for any $\mathbb{E}$-triangle
$$A\stackrel{x}\longrightarrow B\stackrel{y}\longrightarrow C\stackrel{\delta}\dashrightarrow$$ in $\xi$
and any morphism  $a:A \rightarrow A'$, then any $\mathbb{E}$-triangle
$A'\stackrel{x'}\longrightarrow B'\stackrel{y'}\longrightarrow C\stackrel{a_{*}\delta}\dashrightarrow$ belongs to $\xi$.

A class of $\mathbb{E}$-triangles $\xi$ is called $saturated$ if in the situation of Proposition \cite[Proposion 3.15]{NP_2019}, whenever\\
$A_{2}\stackrel{x_{2}}\longrightarrow B_{2}\stackrel{y_{2}}\longrightarrow C\stackrel{\delta_{2}}\dashrightarrow$ and
$A_{1}\stackrel{m_{1}}\longrightarrow M\stackrel{e_{1}}\longrightarrow B_{2}\stackrel{y_{2}^{*}\delta_{1}}\dashrightarrow$ belong to $\xi$, then the $\mathbb{E}$-triangle
$A_{1}\stackrel{x_{1}}\longrightarrow B_{1}\stackrel{y_{1}}\longrightarrow C\stackrel{\delta_{1}}\dashrightarrow$ belongs to $\xi$.

\begin{definition}\cite[Definition 3.1]{JSH1_2020} Let $\xi$ be a class of $\mathbb{E}$-triangles which is closed under isomorphisms. $\xi$ is called a proper class of $\mathbb{E}$-triangles if the following conditions hold:

(1) $\xi$ is closed under finite coproducts and $\Delta_{0} \subseteq \xi$.

(2) $\xi$ is closed under base change and cobase change.

(3) $\xi$ is saturated.
\end{definition}

\begin{definition}\cite[Definition 4.1]{JSH1_2020} An object $P \in \mathcal{C}$ is called $\xi$-$projective$ if for any $\mathbb{E}$-triangle
$$A\stackrel{x}\longrightarrow B\stackrel{y}\longrightarrow C\stackrel{\delta}\dashrightarrow$$ in $\xi$,
the induced sequences of abelian groups
$$0 \longrightarrow \mathcal{C}(P,A)\stackrel{}\longrightarrow \mathcal{C}(P,B)\stackrel{}\longrightarrow \mathcal{C}(P,C)\longrightarrow 0$$ is exact. Dually, we have the definition of $\xi$-$injective$ objects.
\end{definition}

We denote by $\mathcal{P}(\xi)$ (resp. $\mathcal{I}(\xi)$) the class of $\xi$-projective (resp. $\xi$-injective) objects of $\mathcal{C}$. It follows from the definition that the subcategories $\mathcal{P}(\xi)$ and $\mathcal{I}(\xi)$ are full, additive, closed under isomorphisms and direct summands.

An extriangulated category $(\mathcal{C},\mathbb{E},\mathfrak{s})$ is said to have $enough$ $\xi$-$projectives$ (resp. $enough$ $\xi$-$injectives$) provided that for each object $A$ there exists an $\mathbb{E}$-triangle $K\stackrel{}\longrightarrow P\stackrel{}\longrightarrow A\stackrel{}\dashrightarrow$ (resp. $A\stackrel{}\longrightarrow I\stackrel{}\longrightarrow K\stackrel{}\dashrightarrow$ in $\xi$ with $P \in \mathcal{P}(\xi)$ (resp. $I \in \mathcal{I}(\xi)$).

Let $K\stackrel{}\longrightarrow P\stackrel{}\longrightarrow A\stackrel{}\dashrightarrow$ be an $\mathbb{E}$-triangle in $\xi$ with $P \in \mathcal{P}(\xi)$, then we call $K$ the $first$ $\xi$-$syzygy$ of $A$. An $n$th $\xi$-$syzygy$ of $A$ is defined as usual by induction. By Schanuel's lemma (\cite[Proposition 4.3]{JSH1_2020}), any two $\xi$-$syzygies$ of $A$ are isomorphic modulo $\xi$-projectives.

The $\xi$-$projective$ $dimension$ $\xi$-pd$A$ of $A \in \mathcal{C}$ is defined inductively. If $A \in \mathcal{P}(\xi)$, then define $\xi$-pd$A$ = 0. Next if $\xi$-pd$A > 0$, define $\xi$-pd$A \leq n$ if there exists an $\mathbb{E}$-triangle $K\stackrel{}\longrightarrow P\stackrel{}\longrightarrow A\stackrel{}\dashrightarrow$ in $\xi$ with $P \in \mathcal{P}(\xi)$ and $\xi$-pd$K \leq n-1$. Finally we define $\xi$-pd$A = n$ if $\xi$-pd$A \leq n$ and $\xi$-pd$A \nleq n-1$. Of course we set $\xi$-pd$A = \infty$, if $\xi$-pd$A \neq n$ for all $n \geq 0$.

Dually, we can define the $\xi$-$injective$ $dimension$ $\xi$-id$A$ of an object $A \in \mathcal{C}$.

\begin{definition}\cite[Definition 4.4]{JSH1_2020} A $\xi$-$exact$ $complex$ $\mathbf{X}$ is a diagram
$$\cdots \longrightarrow X_{1}\stackrel{d_{1}} \longrightarrow X_{0}\stackrel{d_{0}} \longrightarrow X_{-1} \longrightarrow \cdots$$
in $\mathcal{C}$ such that for each integer $n$, there exists an $\mathbb{E}$-triangle
$K_{n+1}\stackrel{g_{n}}\longrightarrow X_{n}\stackrel{f_{n}}\longrightarrow K_{n}\stackrel{\delta_{n}}\dashrightarrow$ in $\xi$ and $d_{n} =g_{n-1}f_{n}$.
\end{definition}

\begin{definition}\cite[Definition 4.5]{JSH1_2020} Let $\mathcal{W}$ be a class of objects in $\mathcal{C}$. An $\mathbb{E}$-triangle
$$A\stackrel{}\longrightarrow B\stackrel{}\longrightarrow C\stackrel{\delta}\dashrightarrow$$ in $\xi$ is called to be $\mathcal{C}(-, \mathcal{W})-exact$ (resp. $\mathcal{C}(\mathcal{W}, -)-exact$) if for any $W \in \mathcal{W}$, the induced sequences of abelian groups
$0 \longrightarrow \mathcal{C}(C,W)\stackrel{}\longrightarrow \mathcal{C}(B,W)\stackrel{}\longrightarrow \mathcal{C}(A,W)\longrightarrow 0$
(resp. \\ $0 \longrightarrow \mathcal{C}(W,A)\stackrel{}\longrightarrow \mathcal{C}(W,B)\stackrel{}\longrightarrow \mathcal{C}(W,C)\longrightarrow 0$) is exact in $\mathcal{A}b$.
\end{definition}

\begin{definition}\cite[Definition 4.6]{JSH1_2020} Let $\mathcal{W}$ be a class of objects in $\mathcal{C}$. A complex $\mathbf{X}$ is called $\mathcal{C}(-, \mathcal{W})-exact$ (resp. $\mathcal{C}(\mathcal{W}, -)-exact$) if it is a $\xi$-exact complex
$$\cdots \longrightarrow X_{1}\stackrel{d_{1}} \longrightarrow X_{0}\stackrel{d_{0}} \longrightarrow X_{-1} \longrightarrow \cdots$$
in $\mathcal{C}$ such that there is a $\mathcal{C}(-, \mathcal{W})-exact$ (resp. $\mathcal{C}(\mathcal{W}, -)-exact$) $\mathbb{E}$-triangle
$K_{n+1}\stackrel{g_{n}}\longrightarrow X_{n}\stackrel{f_{n}}\longrightarrow K_{n}\stackrel{\delta_{n}}\dashrightarrow$ in $\xi$ for each integer $n$ and $d_{n} =g_{n-1}f_{n}$.
\end{definition}

\begin{definition}\cite[Definition 3.1]{JSH2_2020} Let $M$ be an object in $\mathcal{C}$. A $\xi$-$projective~resolution$ of $M$ is a $\xi$-exact complex $\xymatrix@C=0.5cm{
  \mathbf{P} \ar[r]^{} & M  }$
such that $P_{n} \in \mathcal{P}(\xi)$ for all $n \geq 0$. Dually,  A $\xi$-$injective~coresolution$ of $M$ is a $\xi$-exact complex $\xymatrix@C=0.5cm{
  M \ar[r]^{} & \mathbf{I}  }$
such that $I_{n} \in \mathcal{I}(\xi)$ for all $n \leq 0$.

\end{definition}

\begin{definition}\cite[Definition 3.2]{JSH2_2020} Let $M$ $N$ be objects in $\mathcal{C}$.

(1) If we choose a $\xi$-projective~resolution $\xymatrix@C=0.5cm{
  \mathbf{P} \ar[r]^{} & M  }$ of M, then for any integer $n \geq 0$, the $\xi$-$cohomology~groups$ $\xi xt_{\mathcal{P}(\xi)}^{n}(M,N)$ are defined as
 $$\xi xt_{\mathcal{P}(\xi)}^{n}(M,N) = H^{n}(\mathcal{C}(\mathbf{P},N)).$$

(2) If we choose a $\xi$-injective~coresolution $\xymatrix@C=0.5cm{
   M \ar[r]^{} & \mathbf{I}  }$ of M, then for any integer $n \geq 0$, the $\xi$-$cohomology~groups$ $\xi xt_{\mathcal{I}(\xi)}^{n}(M,N)$ are defined as
 $$\xi xt_{\mathcal{I}(\xi)}^{n}(M,N) = H^{n}(\mathcal{C}(M,\mathcal{I})).$$

\end{definition}

\begin{remark} By \cite[Lemma 3.2]{JSH3_2020}, one can see that $\xi xt_{\mathcal{P}(\xi)}^{n}(-,-)$ and $\xi xt_{\mathcal{I}(\xi)}^{n}(-,-)$ are cohomological functors for any integer $n \geq 0$, independent of the choice of $\xi$-projective~resolutions and $\xi$-injective~resolutions, respectively. In fact, with the modifications of the usual proof, one obtains the isomorphism $\xi xt_{\mathcal{P}(\xi)}^{n}(M,N) \cong \xi xt_{\mathcal{I}(\xi)}^{n}(M,N)$, which is denoted by $\xi xt_{(\xi)}^{n}(M,N)$.
\end{remark}

\section{Main results}

Throughout this section, we assume that $\mathcal{C} = (\mathcal{C},\mathbb{E},\mathfrak{s})$ is an extriangulated category satisfying Condition(WIC) and $\xi$ is a proper class of $\mathbb{E}$-triangles in $\mathcal{C}$. We also assume that $\mathcal{C}$ admits any coproducts, having enough $\xi$-projectives and $\xi$-injectives. In what follows, we also assume that $\mathcal{P}(\xi)$ is a generating subcategory of $\mathcal{C}$.

Let $\mathcal{X}$ and $\mathcal{Y}$ be classes of objects of $\mathcal{C}$. We recall the following right and left orthogonal classes:
$$\mathcal{X}^{\perp} = \{Y \in \mathcal{C} \mid \xi \xt_{\xi}^{1}(X,Y) = 0, \for~ \any~ X \in \mathcal{X}\},$$
$$^{\perp}\mathcal{Y} = \{X \in \mathcal{C} \mid \xi \xt_{\xi}^{1}(X,Y) = 0, \for~ \any~ Y \in \mathcal{Y}\}.$$
We denote by $\Add \mathcal{X}$ the subcategory of all summands of direct sums of objects in $\mathcal{X}$. We write 
$\Pres_{\mathcal{P}(\xi)}^{1}(\Add T) = \{M \in \mathcal{C} \mid$ there exists an $\mathbb{E}$-triangle
$K_{1}\stackrel{}\longrightarrow Y_{0}\stackrel{}\longrightarrow M\stackrel{}\dashrightarrow$ lies in $\xi$ with $Y_{0} \in \Add T$\}. \\

Now, we first give the definition of $\xi$-$tilting~ object$ in an extriangulated category in the context we are working with.

\begin{definition}\label{D31} Let $T$ be a non-zero object in $\mathcal{C}$. $T$ is said to be an $\xi$-$tilting~ object$ if it satisfies the following conditions.

(T1) $\xi$-pdT $\leq 1$.

(T2) $\xi xt_{\xi}^{1}(T,T^{(\lambda)}) = 0$ for any cardinal $\lambda$. In other words, $T$ is a self-orthogonal object in $\mathcal{C}$.

(T3) For any object $P \in \mathcal{P}(\xi)$, there exists an $\mathbb{E}$-triangle $$P\stackrel{}\longrightarrow T_{0}\stackrel{}\longrightarrow T_{1}\stackrel{}\dashrightarrow$$ in $\xi$ with $T_{0}, T_{1} \in \Add T$.

Also, $T$ is called a partial $\xi$-$tilting~ object$ if it satisfies conditions above (T1) and (T2).
\end{definition}

This definition generalizes the classical ones for artin algebras given by Happel and Ringel \cite{DH_1982} and Miyashita \cite{YM_1986}, as well as the one given by Y.G. Hu \cite{YGH_2020} for $\xi$-tilting object in a triangulated category.

The following lemma is frequently used.

\begin{lemma}\label{L32} Let $T$ be an object in $\mathcal{C}$. Then the following statements holds.

(1) If the $\mathbb{E}$-triangle $A\stackrel{}\longrightarrow B\stackrel{}\longrightarrow C\stackrel{}\dashrightarrow$ lies in $\xi$ such that
$B \in \Pres_{\mathcal{P}(\xi)}^{1}(\Add T)$, then $C \in \Pres_{\mathcal{P}(\xi)}^{1}(\Add T)$.

(2) If $\Pres_{\mathcal{P}(\xi)}^{1}(\Add T) = T^{\perp}$, then $T$ is a partial $\xi$-tilting object.
\end{lemma}

\noindent{\textbf{Proof}}
(1) Suppose that the $\mathbb{E}$-triangle
$A\stackrel{x}\longrightarrow B\stackrel{y}\longrightarrow C\stackrel{\delta}\dashrightarrow$ lies in $\xi$ such that
$B \in \Pres_{\mathcal{P}(\xi)}^{1}(\Add T)$. Then there exists an $\mathbb{E}$-triangle
$K_{B}\stackrel{x_{1}}\longrightarrow T_{0}\stackrel{y_{1}}\longrightarrow B\stackrel{\delta_{1}}\dashrightarrow$ in $\xi$ with $T_{0} \in \Add T$.
Since  $\mathcal{C}$ has enough $\xi$-projectives, there exists an $\mathbb{E}$-triangle
$K_{A}\stackrel{x_{2}}\longrightarrow P_{A}\stackrel{y_{2}}\longrightarrow A\stackrel{\delta_{2}}\dashrightarrow$ in $\xi$ with $P_{A} \in \mathcal{P}(\xi)$. We obtain a commutative diagram of $\mathbb{E}$-triangles
$$\xymatrix{
K_{A}\ar[d]^{x_{2}} & K_{B}\oplus P_{A}\ar[d]^{\left[\begin{smallmatrix} x_{1} & 0 \\ 0 & 1 \end{smallmatrix}\right]} & \\
P_{A}\ar[d]^{y_{2}}\ar[r]^{\left[\begin{smallmatrix} 0 \\ 1 \end{smallmatrix}\right]} & T_{0}\oplus P_{A}\ar[r]^{\left[\begin{smallmatrix}   1 & 0 \end{smallmatrix}\right]}\ar[d]^{\left[\begin{smallmatrix} y_{1} & xy_{2} \end{smallmatrix}\right]} & T_{0}\ar@{-->}[r] & \\
A\ar@{-->}[d]\ar[r]_{x} & B\ar@{-->}[d]\ar[r]_{y} & C\ar@{-->}[r]_{\delta} & \\
&&&
}
$$
Since $\mathcal{C}$ satisfies Condition(WIC), by \cite[Lemma 5.9]{NP_2019}, then for some $X \in \mathcal{C}$, we obtain $\mathbb{E}$-triangles \\
$K_{A}\stackrel{}\longrightarrow K_{B} \oplus P_{A}\stackrel{}\longrightarrow X\stackrel{}\dashrightarrow$ and
$X\stackrel{}\longrightarrow T_{0}\stackrel{yy_{1}}\longrightarrow C\stackrel{}\dashrightarrow$.

We now show that $X\stackrel{}\longrightarrow T_{0}\stackrel{yy_{1}}\longrightarrow C\stackrel{}\dashrightarrow$ lies in $\xi$. Since $A\stackrel{x}\longrightarrow B\stackrel{y}\longrightarrow C\stackrel{\delta}\dashrightarrow$ lies in $\xi$ and $K_{B}\stackrel{x_{1}}\longrightarrow T_{0}\stackrel{y_{1}}\longrightarrow B\stackrel{\delta_{1}}\dashrightarrow$ in $\xi$ with $T_{0} \in \Add T$, by \cite[Corollary 3.5]{JSH1_2020}, the class of $\xi$-deflations is closed under compositions, then the $\mathbb{E}$-triangle $X\stackrel{}\longrightarrow T_{0}\stackrel{yy_{1}}\longrightarrow C\stackrel{}\dashrightarrow$ lies in $\xi$. Therefore,
$C \in \Pres_{\mathcal{P}(\xi)}^{1}(\Add T)$ by the definition.

(2) It suffices to show that $\xi$-pdT $\leq 1$. For any object $M \in \mathcal{C}$, since $\mathcal{C}$ has enough $\xi$-injectives, there exists an $\mathbb{E}$-triangle
$M\stackrel{x_{2}}\longrightarrow I\stackrel{y_{2}}\longrightarrow K\stackrel{\delta_{2}}\dashrightarrow$ in $\xi$ with $I \in \mathcal{I}(\xi)$.
Note that, $\mathcal{I}(\xi) \subseteq T^{\perp} = \Pres_{\mathcal{P}(\xi)}^{1}(\Add T)$, then $I \in \Pres_{\mathcal{P}(\xi)}^{1}(\Add T)$. By (1), we have $K \in \Pres_{\mathcal{P}(\xi)}^{1}(\Add T)$. Applying the functor $\mathcal{C}(T,-)$ to the $\mathbb{E}$-triangle
$M\stackrel{x_{2}}\longrightarrow I\stackrel{y_{2}}\longrightarrow K\stackrel{\delta_{2}}\dashrightarrow$, by \cite[Lemma 3.4]{JSH2_2020}, we have an exact sequence
$$0 = \xi xt_{\xi}^{1}(T,K) \stackrel{ }\longrightarrow \xi xt_{\xi}^{2}(T,M)\stackrel{ }\longrightarrow \xi xt_{\xi}^{2}(T,I) = 0.$$
By \cite[Lemma 3.9 (1)]{JSH3_2020}, since $\mathcal{P}(\xi)$ is a generating subcategory of $\mathcal{C}$, then $\xi$-pdT $\leq 1$ if and only if $\xi xt_{\xi}^{2}(T,M) = 0$ for any object $M \in \mathcal{C}$. The proof is completed.

\hfill$\Box$

\begin{lemma}\label{L33} Let $T$ be an object in $\mathcal{C}$. If $\Pres_{\mathcal{P}(\xi)}^{1}(\Add T) = T^{\perp}$, then for each $M \in \Pres_{\mathcal{P}(\xi)}^{1}(\Add T)$, there exists an $\mathbb{E}$-triangle
$K\stackrel{}\longrightarrow T_{M}\stackrel{}\longrightarrow M\stackrel{}\dashrightarrow$ lies in $\xi$
with $T_{M} \in \Add T$ and $K \in \Pres_{\mathcal{P}(\xi)}^{1}(\Add T)$.
\end{lemma}

\noindent{\textbf{Proof}}
For each $M \in \Pres_{\mathcal{P}(\xi)}^{1}(\Add T)$, there exists an $\mathbb{E}$-triangle
$K\stackrel{x}\longrightarrow T_{M}\stackrel{y}\longrightarrow M\stackrel{\delta_{1}}\dashrightarrow$ lies in $\xi$
with $T_{M} \in \Add T$.
Since  $\mathcal{C}$ has enough $\xi$-projectives, there exists an $\mathbb{E}$-triangle
$K_{1}\stackrel{x_{1}}\longrightarrow P_{0}\stackrel{y_{1}}\longrightarrow T\stackrel{\delta_{1}}\dashrightarrow$ in $\xi$ with $P_{0} \in \mathcal{P}(\xi)$. It follows from Lemma \ref{L32} that $T$ is self-orthogonal and $\xi$-pdT $\leq 1$, by the definition of $\xi$-projective dimension, we obtain that $\xi$-pd$K_{1}$ $\leq 0$, thus $K_{1} \in \mathcal{P}(\xi)$.

Applying the functor $\mathcal{C}(P_{0},-)$ to the $\mathbb{E}$-triangle
$K\stackrel{x}\longrightarrow T_{M}\stackrel{y}\longrightarrow M\stackrel{\delta_{1}}\dashrightarrow$, by \cite[Lemma 3.4]{JSH2_2020},
We have the following commutative diagram
$$\xymatrix{
  0 \ar@{-->}[r]& \mathcal{C}(P_{0},K) \ar[d]_{\cong} \ar[r]^{}     & \mathcal{C}(P_{0},T_{M}) \ar[d]_{\cong} \ar[r]^{}     & \mathcal{C}(P_{0},M) \ar[d]_{\cong} \ar@{-->}[r]     & 0\\
  0 \ar[r]^{} & \xi xt_{\xi}^{0}(P_{0},K) \ar[r]^{} & \xi xt_{\xi}^{0}(P_{0},T_{M}) \ar[r]^{} & \xi xt_{\xi}^{0}(P_{0},M) \ar[r]^{}     & \xi xt_{\xi}^{1}(P_{0},K) = 0  }
$$
then the first row in above diagram is exact.

Applying the functor $\mathcal{C}(K_{1},-)$ to the $\mathbb{E}$-triangle
$K\stackrel{x}\longrightarrow T_{M}\stackrel{y}\longrightarrow M\stackrel{\delta_{1}}\dashrightarrow$, by \cite[Lemma 3.4]{JSH2_2020} again,
We have the following commutative diagram
$$\xymatrix{
  0 \ar@{-->}[r]& \mathcal{C}(K_{1},K) \ar[d]_{\cong} \ar[r]^{}     & \mathcal{C}(K_{1},T_{M}) \ar[d]_{\cong} \ar[r]^{}     & \mathcal{C}(K_{1},M) \ar[d]_{\cong} \ar@{-->}[r]     & 0\\
  0 \ar[r]^{} & \xi xt_{\xi}^{0}(K_{1},K) \ar[r]^{} & \xi xt_{\xi}^{0}(K_{1},T_{M}) \ar[r]^{} & \xi xt_{\xi}^{0}(K_{1},M) \ar[r]^{}     & \xi xt_{\xi}^{1}(K_{1},K) = 0  }
$$
then the first row in above diagram is exact.

Since $M \in T^{\perp}$, applying the functor $\mathcal{C}(-,M)$ to the $\mathbb{E}$-triangle
$K_{1}\stackrel{x_{1}}\longrightarrow P_{0}\stackrel{y_{1}}\longrightarrow T\stackrel{\delta_{1}}\dashrightarrow$ , by \cite[Lemma 3.4]{JSH2_2020} again,
We have the following commutative diagram
$$\xymatrix{
    & \mathcal{C}(T,M) \ar[d]_{} \ar[r]^{}     & \mathcal{C}(P_{0},M) \ar[d]_{\cong} \ar[r]^{}     & \mathcal{C}(K_{1},M) \ar[d]_{\cong} \ar@{-->}[r]     & 0\\
  0 \ar[r]^{} & \xi xt_{\xi}^{0}(T,M) \ar[r]^{} & \xi xt_{\xi}^{0}(P_{0},M) \ar[r]^{} & \xi xt_{\xi}^{0}(K_{1},M) \ar[r]^{}     & \xi xt_{\xi}^{1}(T,M) = 0  }
$$
then $\mathcal{C}(P_{0},M) \longrightarrow \mathcal{C}(K_{1},M)$ is an epimorphism.

Since $T$ is self-orthogonal, applying the functor $\mathcal{C}(-,T_{M})$ to the $\mathbb{E}$-triangle
$K_{1}\stackrel{x_{1}}\longrightarrow P_{0}\stackrel{y_{1}}\longrightarrow T\stackrel{\delta_{1}}\dashrightarrow$ , by \cite[Lemma 3.4]{JSH2_2020} again,
We have the following commutative diagram
$$\xymatrix{
    & \mathcal{C}(T,T_{M}) \ar[d]_{} \ar[r]^{}     & \mathcal{C}(P_{0},T_{M}) \ar[d]_{\cong} \ar[r]^{}     & \mathcal{C}(K_{1},T_{M}) \ar[d]_{\cong} \ar@{-->}[r]     & 0\\
  0 \ar[r]^{} & \xi xt_{\xi}^{0}(T,T_{M}) \ar[r]^{} & \xi xt_{\xi}^{0}(P_{0},T_{M}) \ar[r]^{} & \xi xt_{\xi}^{0}(K_{1},T_{M}) \ar[r]^{}     & \xi xt_{\xi}^{1}(T,T_{M}) = 0  }
$$
then $\mathcal{C}(P_{0},T_{M}) \longrightarrow \mathcal{C}(K_{1},T_{M})$ is an epimorphism.

Since the functor $\mathcal{C}(-,-)$ is a biaddtive functor, we have the following commutative diagram
$$\xymatrix{
     & \mathcal{C}(T,K) \ar[d]_{} \ar[r]^{}     & \mathcal{C}(T,T_{M}) \ar[d]_{} \ar[r]^{}     & \mathcal{C}(T,M) \ar[d]_{}        & & \\
  0 \ar[r]^{} &\mathcal{C}(P_{0},K) \ar[d]_{} \ar[r]^{} & \mathcal{C}(P_{0},T_{M}) \ar[d]_{} \ar[r]^{} & \mathcal{C}(P_{0},M) \ar[d]_{} \ar[r]^{}     & 0  &  \\
  0 \ar[r]^{} & \mathcal{C}(K_{1},K) \ar@{-->}[d]_{} \ar[r]^{}     & \mathcal{C}(K_{1},T_{M}) \ar[d]_{} \ar[r]^{} & \mathcal{C}(K_{1},M) \ar[d]_{} \ar[r]^{}     & 0  &  \\
     & 0    & 0  &  0  &  &  }
$$
in which all horizontals and the third and the fourth vertials are exact.

Since $\mathcal{C}(P_{0},T_{M}) \longrightarrow \mathcal{C}(K_{1},T_{M})$ is an epimorphism, for any $g \in \mathcal{C}(K_{1},K)$, there exists $g_{1} \in \mathcal{C}(P_{0},T_{M})$ such that $xg = g_{1}x_{1}$, by (ET3), then we have the following commutative
$$\scalebox{0.9}[1.0]{\xymatrixcolsep{4pc}\xymatrix{
K_{1}\ar[d]_-{g }\ar[r]^{x_{1}}  &          P_{0}\ar@{..>}[ld]_-{h_{1}} \ar[d]^{g_{1} }\ar[r]^-{y_{1}}   & T\ar@{..>}[ld]_-{h_{2}} \ar@{..>}[d]_-{g_{2} } \ar@{-->}[r]^{\delta_{1}} & \\
K\ar[r]^-{x }&          T_{M} \ar[r]^-{ y}   & M \ar@{-->}[r]^{\delta} &.
 }}
$$
Since $T_{M} \in \Add T$, then $g_{2}$ factors through $y$, hence $g$ factors through $x_{1}$ by Lemma \ref {L3}, that is, the morphism  $\mathcal{C}(P_{0},K) \longrightarrow \mathcal{C}(K_{1},K)$ is an epimorphism.

Applying the functor $\mathcal{C}(-,K)$ to the $\mathbb{E}$-triangle
$K_{1}\stackrel{x_{1}}\longrightarrow P_{0}\stackrel{y_{1}}\longrightarrow T\stackrel{\delta_{1}}\dashrightarrow$ , by \cite[Lemma 3.4]{JSH2_2020} again,
We have the following commutative diagram
$$\xymatrix{
    & \mathcal{C}(T,K) \ar[d]_{} \ar[r]^{}     & \mathcal{C}(P_{0},K) \ar[d]_{\cong} \ar[r]^{}     & \mathcal{C}(K_{1},K) \ar[d]_{\cong} \ar[r]    & 0 &\\
  0 \ar[r]^{} & \xi xt_{\xi}^{0}(T,K) \ar[r]^{} & \xi xt_{\xi}^{0}(P_{0},K) \ar[r]^{} & \xi xt_{\xi}^{0}(K_{1},K) \ar[r]^{}     & \xi xt_{\xi}^{1}(T,K) \ar[r]^{} & 0  }
$$
then $\mathcal{\xi} xt_{\xi}^{1}(P_{0},K) \longrightarrow \mathcal{\xi} xt_{\xi}^{1}(K_{1},K)$ is an epimorphism. It implies that $\mathcal{\xi} xt_{\xi}^{1}(T,K) = 0$. Therefore $K \in T^{\perp} = \Pres_{\mathcal{P}(\xi)}^{1}(\Add T)$.

\hfill$\Box$

Now we give an important characterization of the $\xi$-tilting objects: the Bazzoni characterization of $\xi$-tilting object in extriangulated category.

\begin{theorem}\label{T34} Let $T$ be an object in $\mathcal{C}$. Then the following statements are equivalent.

(1) $T$ is an $\xi$-tilting object.

(2) $T^{\perp} = \Pres_{\mathcal{P}(\xi)}^{1}(\Add T)$.
\end{theorem}

\noindent{\textbf{Proof}}
(1) $\Rightarrow$ (2) For each $M \in \Pres_{\mathcal{P}(\xi)}^{1}(\Add T)$, there exists an $\mathbb{E}$-triangle
$K\stackrel{ }\longrightarrow T_{M}\stackrel{ }\longrightarrow M\stackrel{ }\dashrightarrow$ lies in $\xi$
with $T_{M} \in \Add T$. Applying the functor $\mathcal{C}(T,-)$ to the above $\mathbb{E}$-triangle, by \cite[Lemma 3.4]{JSH2_2020}, we have an exact sequence
$$\xi xt_{\xi}^{1}(T,T_{M}) \stackrel{ }\longrightarrow \xi xt_{\xi}^{1}(T,M)\stackrel{ }\longrightarrow \xi xt_{\xi}^{2}(T,K).$$

Since $T$ satisfies the conditions (T1) and (T2), $\xi xt_{\xi}^{1}(T,T_{M}) = \xi xt_{\xi}^{1}(T,K) = 0$. It yields that $\xi xt_{\xi}^{2}(T,M) = 0$ and hence $M \in T^{\perp}$. Therefore, we have that $\Pres_{\mathcal{P}(\xi)}^{1}(\Add T) \subseteq T^{\perp}$.

Now, assume that $M \in T^{\perp}$. Since  $\mathcal{C}$ has enough $\xi$-projectives, there exists an $\mathbb{E}$-triangle
$K\stackrel{x}\longrightarrow P\stackrel{y}\longrightarrow M\stackrel{\delta}\dashrightarrow$ in $\xi$ with $P \in \mathcal{P}(\xi)$.
By (T3), there also exists an $\mathbb{E}$-triangle
$P\stackrel{x_{1}}\longrightarrow T_{0}\stackrel{y_{1}}\longrightarrow T_{1}\stackrel{\delta_{1}}\dashrightarrow$ in $\xi$ with $T_{0}, T_{1} \in \Add T$.
Applying the functor $\mathcal{C}(-,M)$ to the $\mathbb{E}$-triangle
$P\stackrel{x_{1}}\longrightarrow T_{0}\stackrel{y_{1}}\longrightarrow T_{1}\stackrel{\delta_{1}}\dashrightarrow$, by \cite[Lemma 3.4]{JSH2_2020}, we have an exact sequence
$$0 \stackrel{ }\longrightarrow \xi xt_{\xi}^{0}(T_{1},M) \stackrel{ }\longrightarrow \xi xt_{\xi}^{0}(T_{0},M)\stackrel{ }\longrightarrow \xi xt_{\xi}^{0}(P,M) \stackrel{ }\longrightarrow \xi xt_{\xi}^{1}(T_{1},M) = 0$$
By \cite[Lemma 3.4]{JSH2_2020}, We have the following commutative diagram
$$\xymatrix{
  0 \ar@{-->}[r]& \mathcal{C}(T_{1},M) \ar[d]_{ } \ar[r]^{}     & \mathcal{C}(T_{0},M) \ar[d]_{ } \ar[r]^{}     & \mathcal{C}(P,M) \ar[d]_{\cong} \ar@{-->}[r]     & 0\\
  0 \ar[r]^{} & \xi xt_{\xi}^{0}(T_{1},M) \ar[r]^{} & \xi xt_{\xi}^{0}(T_{0},M) \ar[r]^{} & \xi xt_{\xi}^{0}(P,M) \ar[r]^{}     &0 }
$$
Therefore the induced map $\mathcal{C}(T_{0},M) \rightarrow \mathcal{C}(P,M)$ is epic.

Since $\mathcal{C}(T_{0},M) \rightarrow \mathcal{C}(P,M)$ is an epimorphism, there exists $y_{2} \in \mathcal{C}(T_{0},M)$ such that $y = y_{2}x_{1}$. Since $\mathcal{C}$ satisfies Condition (WIC), and $y = y_{2}x_{1}$ is an deflation, then so is $y_{2}$, take $K_{1} =$ CoCone$(y_{2})$. By (ET3)$^{\op}$, we have the following commutative diagram
$$\scalebox{0.9}[1.0]{\xymatrixcolsep{4pc}\xymatrix{
K\ar@{..>}[d]_-{g }\ar[r]^{x}  &          P\ar[d]^{x_{1} }\ar[r]^-{y}   & M \ar@{=}[d]_-{ } \ar@{-->}[r]^{\delta} & \\
K_{1}\ar[r]^-{x_{2}}&          T_{0} \ar[r]^-{y_{2}}   & M \ar@{-->}[r]^{\delta} &.
 }}
$$
Hence $K_{1}\stackrel{x}\longrightarrow T_{0}\stackrel{y}\longrightarrow M\stackrel{\delta}\dashrightarrow$ is an $\mathbb{E}$-triangle in $\xi$ with $T_{0} \in \Add T$ since $\xi$ is closed under cobase change.  Therefore, we have that $T^{\perp} \subseteq \Pres_{\mathcal{P}(\xi)}^{1}(\Add T)$.

(2) $\Rightarrow$ (1) By Lemma \ref{L32}, it suffices to show that $T$ satisfies the condition (T3). For any $P \in P(\xi)$, there exists an $\mathbb{E}$-triangle
$P\stackrel{x}\longrightarrow I\stackrel{y}\longrightarrow K_{0}\stackrel{\delta}\dashrightarrow$ in $\xi$ with $I \in I(\xi)$ since $\mathcal{C}$ has enough $\xi$-injective objects. Note that $I \in I(\xi) \subseteq T^{\perp} = \Pres_{\mathcal{P}(\xi)}^{1}(\Add T)$ and so, there exists an $\mathbb{E}$-triangle
$H\stackrel{x_{1}}\longrightarrow T_{0}\stackrel{y_{1}}\longrightarrow I\stackrel{\delta_{1}}\dashrightarrow$ in $\xi$ with $T_{0} \in \Add T$. By the projectively of $P$, there exists $x_{2}: P \rightarrow T_{0}$ such that $x = y_{1}x_{2}$.
Since $\mathcal{C}$ satisfies Condition (WIC), $x = y_{1}x_{2}$ is an inflation, then so is $x$, take $K_{1} =$ Cone$(x_{2})$. Because $\xi$ is closed under base change, then $P\stackrel{x_{2}}\longrightarrow T_{0}\stackrel{y_{2}}\longrightarrow K_{1}\stackrel{\delta_{2}}\dashrightarrow$ is an $\mathbb{E}$-triangle in $\xi$ with $T_{0} \in \Add T$. This implies that $K_{1} \in \Pres_{\mathcal{P}(\xi)}^{1}(\Add T) = T^{\perp}$.

We obtain a commutative diagram of $\mathbb{E}$-triangles
$$\xymatrix{
0\ar[d]^{ } & H\ar[d]^{x_{1}} & \\
P\ar@{=}[d]^{ }\ar[r]^{x_{2} } & T_{0} \ar[r]^{y_{2} }\ar[d]^{y_{1} } & K_{1}\ar@{-->}[r]^{\delta_{2}} & \\
P\ar@{-->}[d]\ar[r]_{x} & I\ar@{-->}[d]^{\delta_{1}}\ar[r]_{y} & K_{0}\ar@{-->}[r]_{\delta} & .\\
&&&
}
$$
Since $\mathcal{C}$ satisfies Condition(WIC), applying \cite[Lemma 5.9]{NP_2019}, we obtain $\mathbb{E}$-triangles
$0 \stackrel{}\longrightarrow H \stackrel{}\longrightarrow H \stackrel{}\dashrightarrow$ and
$H \stackrel{y_{2}x_{1}}\longrightarrow K_{1}\stackrel{y_{3}}\longrightarrow K_{0}\stackrel{}\dashrightarrow$.

Since $K_{1} \in \Pres_{\mathcal{P}(\xi)}^{1}(\Add T) = T^{\perp}$, there exists an $\mathbb{E}$-triangle
$F\stackrel{x_{3}}\longrightarrow T_{K_{1}}\stackrel{y_{3}}\longrightarrow K_{1}\stackrel{\delta_{3}}\dashrightarrow$ in $\xi$ with $T_{K_{1}} \in \Add T$. According to Lemma \ref{L33}, $F \in \Pres_{\mathcal{P}(\xi)}^{1}(\Add T) = T^{\perp}$.

According to \cite[Proposition 3.15]{NP_2019}, we obtain a commutative diagram of $\mathbb{E}$-triangles
$$\xymatrix{
 & F\ar@{=}[r]^{ }\ar[d]^{ } & F\ar[d]^{x_{3}}  \\
P\ar@{=}[d]^{ }\ar[r]^{ } & E \ar[r]^{ }\ar[d]^{ } & T_{K_{1}}\ar@{-->}[r]^{ }\ar[d]^{y_{3} } & \\
P\ar[r]_{x_{2}} & T_{0}\ar@{-->}[d]^{}\ar[r]_{y_{2}} & K_{1}\ar@{-->}[d]^{\delta_{3}}\ar@{-->}[r]_{\delta_{2}} & .\\
&&&
}
$$
It is easy to see that $\mathbb{E}$-triangles $P\stackrel{}\longrightarrow E \stackrel{ }\longrightarrow T_{K_{1}}\stackrel{}\dashrightarrow$ and $F\stackrel{}\longrightarrow E \stackrel{ }\longrightarrow T_{0}\stackrel{}\dashrightarrow$ in $\xi$, because $\xi$ is closed under base change.

Now we claim that the $\mathbb{E}$-triangle $P\stackrel{}\longrightarrow E \stackrel{ }\longrightarrow T_{K_{1}}\stackrel{}\dashrightarrow$ is the desired $\mathbb{E}$-triangle. Since $F, T_{0} \in T^{\perp}$, then $E \in T^{\perp}$. This implies that there exists an $\mathbb{E}$-triangle $Y\stackrel{}\longrightarrow T_{E} \stackrel{ }\longrightarrow E\stackrel{}\dashrightarrow$ in $\xi$ with $T_{E} \in \Add T$ and $Y \in \Pres_{\mathcal{P}(\xi)}^{1}(\Add T) = T^{\perp}$.
Applying the functor $\mathcal{C}(-,Y)$ to the $\mathbb{E}$-triangle
$P\stackrel{x_{1}}\longrightarrow E\stackrel{y_{1}}\longrightarrow T_{K_{1}}\stackrel{\delta_{1}}\dashrightarrow$, by \cite[Lemma 3.4]{JSH2_2020}, we have an exact sequence
$$0 \stackrel{ }\longrightarrow \xi xt_{\xi}^{1}(T_{K_{1}},Y) \stackrel{ }\longrightarrow \xi xt_{\xi}^{1}(E,Y)\stackrel{ }\longrightarrow \xi xt_{\xi}^{1}(P,Y) = 0.$$
It follows that $\xi xt_{\xi}^{1}(E,Y) = 0$ and so $\mathbb{E}$-triangle $Y\stackrel{}\longrightarrow T_{E} \stackrel{ }\longrightarrow E\stackrel{}\dashrightarrow$ in $\xi$ is splits. Therefore $E \in \Add T$. The proof is completed.

\hfill$\Box$

We write the class
$\Pres_{\mathcal{P}(\xi)}^{2}(\Add T) = \{M \in \mathcal{C} \mid$ there exists an $\mathbb{E}$-triangle
$K_{1}\stackrel{}\longrightarrow Y_{0}\stackrel{}\longrightarrow M\stackrel{}\dashrightarrow$ in $\xi$  \\ with $Y_{0} \in \Add T$ and $K_{1} \in \Pres_{\mathcal{P}(\xi)}^{1}(\Add T) \}.$

\begin{corollary}\label{C35} If $T$ is an $\xi$-tilting object, then $\Pres_{\mathcal{P}(\xi)}^{1}(\Add T)= \Pres_{\mathcal{P}(\xi)}^{2}(\Add T)$.

\end{corollary}
\noindent{\textbf{Proof}} The proof follows directly from definition and Lemma \ref{L33}.

\hfill$\Box$

\begin{definition} Let $\mathcal{X}$ be a subcategory of $\mathcal{C}$. $\mathcal{X}$ is said to be an $\xi$-$covariantly~ finite~ subcategory$ of $\mathcal{C}$, if for any object $M \in \mathcal{C}$, there exists an $\mathbb{E}$-triangle $M\stackrel{}\longrightarrow X_{0} \stackrel{ }\longrightarrow K\stackrel{}\dashrightarrow$ in $\xi$ with $X_{0} \in \mathcal{X}$, such that the induced map $\mathcal{C}(X_{0},X) \rightarrow \mathcal{C}(M,X)$ is an epimorphism for any $X \in \mathcal{X}$.

Moreover, if $K \in ^{\perp}\mathcal{X}$, then $\mathcal{X}$ is called a special $\xi$-$covariantly~ finite~ subcategory$ of $\mathcal{C}$.

\end{definition}

\begin{proposition}\label{P37} If $T$ is a $\xi$-tilting object in $\mathcal{C}$, let $\mathcal{X} = T^{\perp}$, then the following results hold.

(1) $\mathcal{X} = T^{\perp}$ is a special $\xi$-covariantly finite subcategory of $\mathcal{C}$.

(2) For each $K \in ^{\perp}\mathcal{X}$, $\xi$-pd$K \leq 1$.

(3) If the $\mathbb{E}$-triangle $X\stackrel{}\longrightarrow Y \stackrel{ }\longrightarrow Z\stackrel{}\dashrightarrow$ in $\xi$  with $X, Y \in \mathcal{X}$, then $Z\in \mathcal{X}$.

\end{proposition}

\noindent{\textbf{Proof}} (1) For any $M \in \mathcal{C}$, since $\mathcal{C}$ has enough $\xi$-injective objects, there exists an $\mathbb{E}$-triangle
$M\stackrel{}\longrightarrow I\stackrel{}\longrightarrow K\stackrel{}\dashrightarrow$ in $\xi$ with $I \in I(\xi)$. Note that $T$ is a $\xi$-tilting object in $\mathcal{C}$ and $I \in I(\xi) \subseteq T^{\perp} = \mathcal{X}$, according to Lemma \ref{L32} (1), we have $K \in \Pres_{\mathcal{P}(\xi)}^{1}(\Add T) = T^{\perp}$. Since $\mathcal{C}$ has enough $\xi$-projective objects, there exists an $\mathbb{E}$-triangle
$F\stackrel{x_{1}}\longrightarrow P_{0}\stackrel{y_{1}}\longrightarrow K\stackrel{\delta_{1}}\dashrightarrow$ in $\xi$ with
$P_{0} \in \mathcal{P}(\xi)$.
According to \cite[Proposition 3.15]{NP_2019}, we obtain a commutative diagram of $\mathbb{E}$-triangles
$$\xymatrix{
 & F\ar@{=}[r]^{ }\ar[d]^{ } & F\ar[d]^{x_{1}}  \\
M\ar@{=}[d]^{ }\ar[r]^{ } & E \ar[r]^{ }\ar[d]^{ } & P_{0}\ar@{-->}[r]^{ }\ar[d]^{y_{1} } & \\
M\ar[r]_{} & I\ar@{-->}[d]^{}\ar[r]_{} &K\ar@{-->}[d]^{\delta_{1}}\ar@{-->}[r]_{ } & .\\
&&&
}
$$
Then the $\mathbb{E}$-triangles
$M\stackrel{x}\longrightarrow E\stackrel{y}\longrightarrow P_{0}\stackrel{\delta}\dashrightarrow$ and
$F\stackrel{ }\longrightarrow E\stackrel{ }\longrightarrow I\stackrel{ }\dashrightarrow$ in $\xi$, because $\xi$ is closed under base change.

Since $T$ is a $\xi$-tilting object in $\mathcal{C}$, there exists an $\mathbb{E}$-triangle
$P_{0}\stackrel{x_{2}}\longrightarrow T_{0}\stackrel{y_{2}}\longrightarrow T_{1}\stackrel{\delta_{2}}\dashrightarrow$ in $\xi$ with $T_{0}, T_{1} \in \Add T$. According to (ET4), we obtain a commutative diagram of $\mathbb{E}$-triangles
$$\xymatrix{
F\ar@{=}[d]^{ }\ar[r]^{x_{1} }  & P_{0}\ar[r]^{y_{1}}\ar[d]^{x_{2} } & K\ar[d]^{ }\ar@{-->}[r]^{\delta_{1}} & \\
F\ar[r]^{x_{2}x_{1} } & T_{0} \ar[r]^{y_{3} }\ar[d]^{y_{2} } & G\ar@{-->}[r]^{\delta_{3} }\ar[d]^{  } & \\
  & T_{1}\ar@{-->}[d]^{\delta_{2}}\ar@{=}[r]_{ } &T_{1}\ar@{-->}[d]^{} & .\\
&&& }
$$
According to \cite[Corollary 3.5]{JSH1_2020}, the class of $\xi$-inflations is closed under compositions, then the $\mathbb{E}$-triangle
$F\stackrel{x_{2}x_{1}}\longrightarrow T_{0}\stackrel{y_{3}}\longrightarrow G\stackrel{\delta_{3}}\dashrightarrow$ lies in $\xi$.
Note that $T_{0} \in \Add T$, then $G \in \Pres_{\mathcal{P}(\xi)}^{1}(\Add T) = T^{\perp}$. According to Lemma \ref{L33}, $F \in \Pres_{\mathcal{P}(\xi)}^{1}(\Add T) = T^{\perp}$.
Consider the $\mathbb{E}$-triangle $F\stackrel{x_{1}}\longrightarrow E\stackrel{y_{1}}\longrightarrow I\stackrel{\delta_{1}}\dashrightarrow$, we have $E \in \mathcal{X} = T^{\perp}$, because $F, I \in T^{\perp}$.

By the preceding discussion, we obtain the $\mathbb{E}$-triangle
$M\stackrel{x}\longrightarrow E\stackrel{y}\longrightarrow P_{0}\stackrel{\delta}\dashrightarrow$ in $\xi$ with $E \in \mathcal{X} = T^{\perp}$. Clearly, by \cite[Lemma 3.9]{JSH3_2020}, $P_{0} \in ^{\perp}\mathcal{X}$. For any $X \in \mathcal{X}$, since $\mathcal{C}$ has enough $\xi$-injective objects, there exists an $\mathbb{E}$-triangle
$X\stackrel{x_{0}}\longrightarrow I_{0}\stackrel{y_{0}}\longrightarrow K_{X}\stackrel{\delta_{0}}\dashrightarrow$ in $\xi$ with $I_{0} \in I(\xi)$.

Since $P_{0} \in \mathcal{P}(\xi)$, we have an exact sequence of abelian groups
$$0 \stackrel{ }\longrightarrow \mathcal{C}(P_{0},X)  \stackrel{ }\longrightarrow \mathcal{C}(P_{0},I_{0})  \stackrel{ }\longrightarrow  \mathcal{C}(P_{0},K_{X}) \stackrel{ }\longrightarrow   0.$$

Since $I_{0} \in \mathcal{I}(\xi)$, we have an exact sequence of abelian groups
$$0 \stackrel{ }\longrightarrow \mathcal{C}(P_{0},I_{0})  \stackrel{ }\longrightarrow \mathcal{C}(E,I_{0})  \stackrel{ }\longrightarrow  \mathcal{C}(M,I_{0},) \stackrel{ }\longrightarrow   0.$$

Because the functor $\mathcal{C}(-,-)$ is a biaddtive functor, we have the following commutative diagram
$$\xymatrix{
 & & 0\ar[d]_{} & & \\
0 \ar[r]^{} & \mathcal{C}(P_{0},X) \ar[d]_{} \ar[r]^{}   & \mathcal{C}(P_{0},I_{0}) \ar[d]_{} \ar[r]^{}    & \mathcal{C}(P_{0},K_{X}) \ar[d]_{}\ar[r]^{}     & 0\\
  & \mathcal{C}(E,X) \ar[d]_{} \ar[r]^{} & \mathcal{C}(E,I_{0}) \ar[d]_{} \ar[r]^{} & \mathcal{C}(E,K_{X}) \ar[d]_{}  &   \\
  & \mathcal{C}(M,X) \ar[r]^{}     & \mathcal{C}(M,I_{0}) \ar[d]_{} \ar[r]^{} & \mathcal{C}(M,K_{X})  & \\
     &    & 0  &    &   }
$$
in which all horizontals and vertials are exact.

Since $\mathcal{C}(E,I_{0}) \longrightarrow \mathcal{C}(M,I_{0})$ is an epimorphism, for any $g \in \mathcal{C}(M,X)$, there exists $g_{1} \in \mathcal{C}(E,I_{0})$ such that $x_{0}g = g_{1}x$, by (ET3), then we have the following commutative
$$\scalebox{0.9}[1.0]{\xymatrixcolsep{4pc}\xymatrix{
M\ar[d]_-{g }\ar[r]^{x}  &          E\ar@{..>}[ld]_-{h_{1}} \ar[d]^{g_{1} }\ar[r]^-{y}   & P_{0}\ar@{..>}[ld]_-{h_{2}} \ar@{..>}[d]_-{g_{2} } \ar@{-->}[r]^{\delta} & \\
X\ar[r]^-{x_{0} }&          I_{0} \ar[r]^-{ y_{0}}   & K_{X} \ar@{-->}[r]^{\delta_{0}} &.
 }}
$$
Since $P_{0} \in \mathcal{P}(\xi)$, then $g_{2}$ factors through $y_{0}$, hence $g$ factors through $x$ by Lemma \ref {L3}, that is, the morphism  $\mathcal{C}(E,X) \longrightarrow \mathcal{C}(M,X)$ is an epimorphism. Therefore, $\mathcal{X} = T^{\perp}$ is a special $\xi$-covariantly finite subcategory of $\mathcal{C}$.

(2) Suppose $Y \in ^{\perp} \mathcal{X}$. For any $M \in \mathcal{C}$, since $\mathcal{C}$ has enough $\xi$-injective objects, there exists an $\mathbb{E}$-triangle
$M\stackrel{}\longrightarrow I\stackrel{}\longrightarrow K\stackrel{}\dashrightarrow$ in $\xi$ with $I \in I(\xi)$.
Applying the functor $\mathcal{C}(Y,-)$ to this $\mathbb{E}$-triangle, by \cite[Lemma 3.4]{JSH2_2020}, we have an exact sequence
$$ \xi xt_{\xi}^{1}(Y,K) \stackrel{ }\longrightarrow \xi xt_{\xi}^{2}(Y,M)\stackrel{ }\longrightarrow \xi xt_{\xi}^{2}(Y,I) = 0.$$
Since $K \in T^{\perp} = \mathcal{X}$, we have $\xi xt_{\xi}^{2}(Y,M) = 0$ and so $\xi$-pd$K \leq 1$.

(3) This result is trivial.

\hfill$\Box$

\begin{lemma}\label{L38} Let $T$ be an object in $\mathcal{C}$ and $\mathcal{X} \subseteq T^{\perp}$. Assume that $T$ satisfies conditions (T1) and (T2). If $\mathcal{X} \subseteq \Pres_{\mathcal{P}(\xi)}^{1}(\Add T)$, then for each $X \in \mathcal{X}$, there exists an $\mathbb{E}$-triangle $K_{X}\stackrel{}\longrightarrow T_{X} \stackrel{ }\longrightarrow X\stackrel{}\dashrightarrow$ in $\xi$  with $T_{X} \in \Add T$ and $K_{X} \in T^{\perp}$.

\end{lemma}

\noindent{\textbf{Proof}} Assume that  $X \in \mathcal{X} \subseteq \Pres_{\mathcal{P}(\xi)}^{1}(\Add T)$, there exists an $\mathbb{E}$-triangle $K_{X}\stackrel{x}\longrightarrow T_{X} \stackrel{y }\longrightarrow X\stackrel{\delta}\dashrightarrow$ in $\xi$  with $T_{X} \in \Add T$.
Since $\mathcal{C}$ has enough $\xi$-projective objects, there exists an $\mathbb{E}$-triangle
$K_{T}\stackrel{x_{1}}\longrightarrow P_{T}\stackrel{y_{1}}\longrightarrow T\stackrel{\delta_{1}}\dashrightarrow$ in $\xi$ with
$P_{T} \in \mathcal{P}(\xi)$. Since $T$ satisfies condition (T1), $\xi$-pdT $\leq 1$, by the definition of $\xi$-projective dimension, we obtain that $\xi$-pd$K_{T}$ $\leq 0$, thus $K_{T} \in \mathcal{P}(\xi)$.

Since $P_{T}, K_{T} \in \mathcal{P}(\xi)$, we have exact sequences of abelian groups
$$0 \stackrel{ }\longrightarrow \mathcal{C}(P_{T},K_{X})  \stackrel{ }\longrightarrow \mathcal{C}(P_{0},T_{X})  \stackrel{ }\longrightarrow  \mathcal{C}(P_{T},X) \stackrel{ }\longrightarrow   0,$$
and
$$0 \stackrel{ }\longrightarrow \mathcal{C}(K_{T},K_{X})  \stackrel{ }\longrightarrow \mathcal{C}(K_{T},T_{X})  \stackrel{ }\longrightarrow  \mathcal{C}(K_{T},X) \stackrel{ }\longrightarrow   0.$$

Since $T$ satisfies condition (T2), applying the functor $\mathcal{C}(-,T_{X})$ to the $\mathbb{E}$-triangle
$K_{T}\stackrel{x_{1}}\longrightarrow P_{T}\stackrel{y_{1}}\longrightarrow T\stackrel{\delta_{1}}\dashrightarrow$ , by \cite[Lemma 3.4]{JSH2_2020},
We have the following commutative diagram
$$\xymatrix{
    & \mathcal{C}(T,T_{X}) \ar[d]_{} \ar[r]^{}     & \mathcal{C}(P_{T},T_{X}) \ar[d]_{\cong} \ar[r]^{}     & \mathcal{C}(K_{T},T_{X}) \ar[d]_{\cong} \ar@{-->}[r]     & 0\\
  0 \ar[r]^{} & \xi xt_{\xi}^{0}(T,T_{X}) \ar[r]^{} & \xi xt_{\xi}^{0}(P_{T},T_{X}) \ar[r]^{} & \xi xt_{\xi}^{0}(K_{T},T_{X}) \ar[r]^{}     & \xi xt_{\xi}^{1}(T,T_{M}) = 0  }
$$
then $\mathcal{C}(P_{T},T_{X}) \longrightarrow \mathcal{C}(K_{T},T_{X})$ is an epimorphism.

Because the functor $\mathcal{C}(-,-)$ is a biaddtive functor, we have the following commutative diagram
$$\xymatrix{
& \mathcal{C}(T,K_{X}) \ar[d]_{} \ar[r]^{}   & \mathcal{C}(T,T_{X}) \ar[d]_{} \ar[r]^{}    & \mathcal{C}(T,X) \ar[d]_{}  & \\
0 \ar[r]^{}  & \mathcal{C}(P_{T},K_{X}) \ar[d]_{} \ar[r]^{} & \mathcal{C}(P_{T},T_{X}) \ar[d]_{} \ar[r]^{} & \mathcal{C}(P_{T},X) \ar[d]_{} \ar[r]^{}  & 0   \\
0 \ar[r]^{}  & \mathcal{C}(K_{T},K_{X}) \ar[r]^{}     & \mathcal{C}(K_{T},T_{X}) \ar[d]_{} \ar[r]^{} & \mathcal{C}(K_{T},X)  \ar[r]^{}  & 0 \\
     &    & 0  &    &   }
$$
in which all horizontals and vertials are exact.

Since $\mathcal{C}(P_{T},T_{X}) \longrightarrow \mathcal{C}(K_{T},T_{X})$  is an epimorphism, for any $g \in \mathcal{C}(K_{T},K_{X})$, there exists $g_{1} \in \mathcal{C}(P_{T},T_{X})$ such that $xg = g_{1}x_{1}$, by (ET3), then we have the following commutative diagram
$$\scalebox{0.9}[1.0]{\xymatrixcolsep{4pc}\xymatrix{
K_{T}\ar[d]_-{g }\ar[r]^{x_{1}}  &          P_{T}\ar@{..>}[ld]_-{h_{1}} \ar[d]^{g_{1} }\ar[r]^-{y_{1}}   & T\ar@{..>}[ld]_-{h_{2}} \ar@{..>}[d]_-{g_{2} } \ar@{-->}[r]^{\delta_{1}} & \\
K_{X}\ar[r]^-{x }  &           T_{X} \ar[r]^-{ y}   &  X \ar@{-->}[r]^{\delta} & .
 }}
$$

Since $T_{X} \in \Add T$, then $g_{2}$ factors through $y$, hence $g$ factors through $x_{1}$ by Lemma \ref {L3}, that is, the morphism  $\mathcal{C}(P_{T},K_{X}) \longrightarrow \mathcal{C}(K_{T},K_{X})$ is an epimorphism.

Applying the functor $\mathcal{C}(-,K_{X})$ to the $\mathbb{E}$-triangle
$K_{1}\stackrel{x_{1}}\longrightarrow P_{0}\stackrel{y_{1}}\longrightarrow T\stackrel{\delta_{1}}\dashrightarrow$, Since $P_{T}, K_{T} \in \mathcal{P}(\xi)$, by \cite[Lemma 3.4]{JSH2_2020} again,
We have the following commutative diagram
$$\xymatrix{
    & \mathcal{C}(T,K_{X}) \ar[d]_{} \ar[r]^{}     & \mathcal{C}(P_{T},K_{X}) \ar[d]_{\cong} \ar[r]^{}     & \mathcal{C}(K_{T},K_{X}) \ar[d]_{\cong} \ar[r]    & 0 &\\
  0 \ar[r]^{} & \xi xt_{\xi}^{0}(T,K_{X}) \ar[r]^{} & \xi xt_{\xi}^{0}(P_{T},K_{X}) \ar[r]^{} & \xi xt_{\xi}^{0}(K_{T},K_{X}) \ar[r]^{}     & \xi xt_{\xi}^{1}(T,K_{X}) \ar[r]^{} & 0  }
$$
then $\mathcal{\xi} xt_{\xi}^{1}(P_{T},K_{X}) \longrightarrow \mathcal{\xi} xt_{\xi}^{1}(K_{T},K_{X})$ is an epimorphism. It implies that $\mathcal{\xi} xt_{\xi}^{1}(T,K_{X}) = 0$. Therefore $K_{X} \in T^{\perp}$.

\hfill$\Box$

\begin{lemma}\label{L39} Let $\mathcal{X} \subseteq \mathcal{C}$ be a class of objects such that $\mathcal{X} \cap ^{\perp}\mathcal{X}$ is closed under coproducts. Suppose $\mathcal{X}$ satisfies the condition that if the $\mathbb{E}$-triangle
$X\stackrel{}\longrightarrow Y \stackrel{ }\longrightarrow Z\stackrel{}\dashrightarrow$ in $\xi$ with $X, Y \in \mathcal{X}$, then $Z \in \mathcal{X}$. Then $\mathcal{X} \cap ^{\perp}\mathcal{X}$ is closed under direct summands.

\end{lemma}

\noindent{\textbf{Proof}}
Using the Eilenberg's swindle \cite[Proposition 1.4]{HH_2004}, one can prove that $\mathcal{X} \cap ^{\perp}\mathcal{X}$ is closed under direct summands. The proof is similar to \cite[Lemma 3.12]{YGH_2020}.

\hfill$\Box$

\begin{proposition}\label{P310} Let $\mathcal{X} \subseteq \mathcal{C}$ be a class of objects such that $\mathcal{X} \cap ^{\perp}\mathcal{X}$ is closed under coproducts. If $\mathcal{X}$ satisfies the following conditions,

(1) $\mathcal{X}$ is a special $\xi$-covariantly finite subcategory in $\mathcal{C}$;

(2) for each $K \in ^{\perp}\mathcal{X}$, $\xi$-pd$K \leq 1$;

(3) If the $\mathbb{E}$-triangle $X\stackrel{}\longrightarrow Y \stackrel{ }\longrightarrow Z\stackrel{}\dashrightarrow$ in $\xi$  with $X, Y \in \mathcal{X}$, then $Z\in \mathcal{X}$. \\
Then there is a $\xi$-tilting object $T$ in $\mathcal{C}$ such that $\mathcal{X} = T^{\perp}$.
\end{proposition}

\noindent{\textbf{Proof}}
The proof is similar to \cite[Lemma 3.13]{YGH_2020}.

\hfill$\Box$


\end{document}